\documentclass{article}
\usepackage{amsfonts}
\usepackage{amssymb}
\usepackage{amsbsy}
\usepackage{amsthm}
\usepackage{amsmath} 
\usepackage{comment}

\newtheorem{theorem}{Theorem}
\newtheorem{lemma}[theorem]{Lemma}

\begin{document}

\title{A generalisation of a congruence due to Vantieghem only holding for primes}
\author{L J P Kilford\footnote{Address: Mathematics 253-37, California Institute of Technology, Pasadena CA 91125.}}

\maketitle
\begin{abstract}
In this note we present a family of congruences which hold if and only if a natural number~$n$ is prime.
\end{abstract}

The subject of primality testing has been in the mathematical and general news recently, with the announcement~\cite{AKS2002} that there exists a polynomial-time algorithm to determine whether an integer~$p$ is prime or not.

There are older deterministic primality tests which are less efficient; the classical example is Wilson's Theorem, that
\[
(n-1)! \equiv -1\mod n\text{ if and only if }n\text{ is prime.}
\]
Although this is a deterministic algorithm, it does not provide a workable primality test because it requires much more calculation than trial division.

This note provides another congruence satisfied by primes and only by primes; it is a generalisation of previous work. In Guy~\cite{guy-unsolved}, problem~$\boldsymbol{A17}$, the following result due to Vantieghem~\cite{vantieghem} is quoted:
\begin{theorem}[Vantieghem, \cite{vantieghem}]
Let~$n$ be a natural number greater than~$2$. Then~$n$ is prime if and only if
\[
\prod_{d=1}^{n-1}(1-2^d) \equiv n \mod (2^n-1).
\]
\end{theorem}

In this note, we will generalise this result to obtain the following theorem:
\begin{theorem}
Let~$m$ and $n$ be natural numbers greater than~$2$.
Then~$n$ is prime if and only if
\begin{equation*}
\sharp:\;\prod_{d=1}^{n-1}(1-m^d) \equiv n \mod \frac{m^n-1}{m-1}.
\end{equation*}
\end{theorem}
\begin{proof}
We follow the method of Vantieghem, using a congruence satisfied by cyclotomic polynomials.
\begin{lemma}[Vantieghem]
\label{vantieghem-lemma}
Let~$m$ be a natural number greater than~$1$ and let~$\Phi_m(X)$ be the~$m^{th}$ cyclotomic polynomial. Then
\[
\underset{(d,m)=1}{\prod_{d=1}^{m}}(X-Y^d)\equiv \Phi_m(X) \mod \Phi_m(Y)\text{ in }\mathbf{Z}[X,Y].
\]
\end{lemma}
\begin{proof}[Proof of Lemma~\ref{vantieghem-lemma}]
We can write
\[
\underset{(d,m)=1}{\prod_{d=1}^{m}}(X-Y^d)-\Phi_m(X)=f_0(Y)+f_1(Y)X+f_2(Y)X^2+\cdots
\]
where the~$f_i$ are polynomials over~$\mathbf{Z}$.

Let~$\zeta$ be a primitive~$m^{th}$ root of unity. Now, if~$Y=\zeta$ then we see that the left hand side of this expression is identically~$0$ in~$X$.

This implies that the~$f_i$ are zero at every~$\zeta$ and every~$i$. Therefore, we have~$f_i(Y) \equiv 0 \mod \Phi_m(Y)$, which is enough to prove the Lemma.
\end{proof}

If~$p$ is prime, then we have that~$\Phi_p(X)=X^{p-1}+X^{p-2}+\cdots+X+1$. Therefore, if we set~$m=p$ in the Lemma, we find that
\[
\prod_{d=1}^{p-1}(X-Y^d)\equiv X^{p-1}+X^{p-2}+\cdots+X+1 \mod (Y^{p-1}+\cdots+1).
\]
We now set~$X=1$ and~$Y=m$, to get
\[
\prod_{d=1}^{p-1}(1-m^d)\equiv p \mod \frac{m^p-1}{m-1};
\]
this proves that if~$p$ is prime then the congruence holds.

We now prove the converse, by supposing that the congruence~$\sharp$ holds, and that~$p$ is not prime. Therefore~$p$ is composite, and hence has a smallest prime factor~$q$. We write~$p=q\cdot a$; now~$q \leq a$, and also~$p \leq a^2$.

Now we have that~$m^a-1$ divides~$m^p-1$ and~$m^a-1$ divides the product~$\prod_{d=1}^{p-1}(m^d-1)$. By combining this with the congruence~$\sharp$ in the Theorem, this implies that~$(m^a-1)/(m-1)$ divides~$p$. Therefore we have
\[
\frac{m^a-1}{m-1} \leq p \leq a^2.
\]
Now this is only possible, when~$m\geq 3$, for~$m=3$ and~$a=2$. It can be easily checked that the congruence does not hold in this case, so we have proved the Theorem.
\end{proof}
Guy also asks if there is a relationship between the congruence given by Vantieghem and Wilson's Theorem. The following theorem gives an elementary congruence similar to that of Vantieghem between a product over integers and a cyclotomic polynomial. It is in fact equivalent to Wilson's Theorem.
\begin{theorem}
\label{wilson-equiv}
Let~$m$ be a natural number greater than~2. Then~$m$ is prime if and only if
\[
\lozenge:\;\Phi_m(X) \equiv F(X):=\underset{(i,m)=1}{\prod_{i=1}^{m-1}} (X-i-1) + 1 \mod m.
\]
\end{theorem}
\begin{proof}[Proof of Theorem~\ref{wilson-equiv}]
Firstly, we prove that if~$m$ is not prime, the congruence~$\lozenge$ in Theorem~\ref{wilson-equiv} does not hold. 

Recall that~$\phi(m)$ is defined to be Euler's totient function; the number of integers in the set~$\left\{1,\ldots,m\right\}$ which are coprime to~$m$.

The coefficient of~$X^{\phi(m)-1}$ on the right-hand side is given by the sum
\[
-\underset{(i,m)=1}{\sum_{i=1}^{m-1}} i+1 = -\phi(m)-\underset{(i,m)=1}{\sum_{i=1}^{m-1}} i \equiv -\phi(m) \mod m;
\]
the final inequality holds because if~$(i,m)=1$ then~$(-i,m)=1$ as well, and the case~$i\equiv -i \mod m$ does not occur because then we have~$2i \equiv 0 \mod m$ and therefore~$2 \equiv 0 \mod m$ which is false because~$m$ is greater than~2.

We now use some theorems to be found in a paper by Gallot~\cite{gallot} (Theorem~1.1 and Theorem~1.4):
\begin{theorem}
\label{gallot-theorem}
Let~$p$ be a prime and~$m$ be a natural number.
\begin{enumerate}
\item The following relations between cyclotomic polynomials hold:
\begin{align*}
\Phi_{pm}(x)&=\Phi_{m}(x^p)\text{ if }p\mid m \\
\Phi_{pm}(x)&=\frac{\Phi_{m}(x^p)}{\Phi_m(x)}\text{ if }p\nmid m.
\end{align*}
\item If~$m > 1$ then
\begin{align*}
\Phi_n(1)&=p\text{ if~$n$ is a power of a prime~$p$}\\
\Phi_n(1)&=1\text{ otherwise.}
\end{align*}
\end{enumerate}
\end{theorem}

From these results, we see that if~$m$ is not a prime power then we have~$\Phi_n(1) \equiv 1 \mod m$, and the right hand side of the congruence~$\lozenge$ when evaluated at~$X=1$ is
\[
1+\underset{(i,m)=1}{\prod_{i=1}^{m-1}} -i 
.
\]
We see that this is not congruent to~$1 \mod m$ because the product is over those~$i$ which are coprime to~$m$, so the product does not vanish modulo~$m$.

If~$m$ is a prime power~$p^n$, then we see from Theorem~\ref{gallot-theorem}.1 that~$\Phi_{p^n}(x)=\Phi_p(x^{p^{n-1}})$; in particular, we see that the coefficient of~$x^{\phi(p^n)-1}$ is~0, which differs from the coefficient of~$x^{\phi(p^n)-1}$ in~$F(X)$. 

Therefore, if~$m$ is not prime then the congruence does not hold. We now show that if~$m$ is prime, the congruence holds.

If~$m$ is prime then~$\Phi_m(x)=x^{m-1}+x^{m-2}+\cdots+x+1$. Let us consider the polynomials~$\Phi_m(X+1)$ and~$F(X+1)$. Now, modulo~$m$ we have
\[
\Phi_m(X+1)=X^{m-1}\text{ and }F(X+1)=\underset{(i,m)=1}{\prod_{i=1}^{m-1}} (X-i) + 1.
\]
Now if~$x \ne 0\mod m$, then we see that~$\Phi_m(x)\equiv 1$ and that~$F(x+1) \equiv 1$, because the product vanishes.

And if we have~$x=0$, then~$\Phi_m(x)=0$ and, by Wilson's Theorem, $F(0) \equiv (m-1)! +1 \equiv 0 \mod m$.

Therefore we have proved the Theorem.
\end{proof}
\providecommand{\bysame}{\leavevmode\hbox to3em{\hrulefill}\thinspace}

\end{document}